\documentclass{article}
\usepackage{amsmath,amssymb,amsthm,amscd, mathtools, latexsym}
\usepackage{enumerate,varioref, dsfont, fancyhdr}
\usepackage{tikz-cd}
\usepackage{enumitem}
\usepackage{multicol}
\usepackage{hyperref}
\usepackage{url}
\usepackage{textcomp}
\usepackage[left=3.5cm, right=3.5cm]{geometry}

\newtheorem{Thm}{Theorem}[section]
\newtheorem{Prop}[Thm]{Proposition}
\newtheorem{Lem}[Thm]{Lemma}
\newtheorem{Cor}[Thm]{Corollary}

\theoremstyle{definition} 
\newtheorem{Def}[Thm]{Definition}
\newtheorem{Rem}[Thm]{Remark}
\newtheorem{Not}[Thm]{Notation}

\newcommand\restr[2]{\ensuremath{\left.#1\right|_{#2}}}

\def\cit{{\mathbb C}}

\def\qit{{\mathbb Q}}
\def\zit{{\mathbb Z}}

\def\git{{\mathbb G}}

\def\0{{\mathcal O}}

\def\Hom{\mathop{\rm Hom}\nolimits}

\def\A{{\mathcal A}}

\begin{document}
\title{Galois descent for higher Brauer groups}

\author{H. Anthony Diaz}
\newcommand{\Addresses}{{\bigskip \footnotesize
\textsc{Department of Mathematics, Washington University, St. Louis, MO 63130} \par \nopagebreak
\textit{Email address}: \ \texttt{humbertoadiaziii@gmail.com}}}

\date{}
\maketitle

\begin{abstract}{\noindent For $X$ a smooth projective variety over a field $k$, we consider the problem of Galois descent for higher Brauer groups. More precisely, we extend a finiteness result of Colliot-Th\'el\`ene and Skorobogatov \cite{CTS} to higher Brauer groups.} 
\end{abstract}

\vspace{5mm}

\noindent For $X$ a smooth projective variety over a field $k$, the Brauer group $Br(X) = H^{2}_{\text{\'et}} (X, \git_{m})$ is a fundamental invariant in arithmetic geometry. Of particular interest is its role in the Tate conjecture in codimension $1$. Recall that the Tate conjecture in codimension $m$, which we denote by $TC^{m} (\overline{X})_{\qit_{\ell}}$, states that when $k$ is a finitely generated field, $\overline{k}$ is its separable closure and $\overline{X} = X \times_{k} \overline{k}$, the cycle class map 
\begin{equation} CH^{m} (\overline{X}) \otimes \qit_{\ell} \to H^{2m}_{\text{\'et}} (\overline{X}, \qit_{\ell}(m))\label{cycle-q} \end{equation}
surjects onto the subspace of {\em Tate classes}:
\[\bigcup_{U} H^{2m}_{\text{\'et}} (\overline{X}, \qit_{\ell}(m))^{U}\]
(where $U$ ranges over all open subgroups of $Gal(\overline{k}/k)$). When $k$ is a finite field, Tate showed in \cite{T} that the Tate conjecture $TC^{1} (\overline{X})_{\qit_{\ell}}$ holds $\Leftrightarrow$ the $\ell$-primary torsion in $Br (X)$ is finite (for $\ell \neq \text{char } k$). For arbitrary fields, the Tate conjecture for divisors is equivalent to the finiteness of the $\ell$-primary torsion in $Br (\overline{X})^{Gal(\overline{k}/k)}$ (see, for instance, \cite{CC} Prop. 2.1.1).\\
\indent In higher codimension, the $m^{th}$ higher Brauer groups $Br^{m}(X)$ are defined by $H^{2m+1}_{L} (X, \zit(m))$, where $H^{*}_{L} (X, \zit(m))$ denote the \'etale motivic cohomology groups. These latter are (\'etale) hyper-cohomology groups of the \'etale sheafification of Bloch’s cycle complexes \cite{Bl}, denoted by $\zit(m)$. When $m=1$, the complex $\zit(1)$ is quasi-isomorphic to $\git_{m}[-1]$, which recovers the usual Brauer group. As motivation for the utility of higher Brauer groups, we begin with the following observation: 
\begin{Prop}\label{Tate} Let $k$ be a finitely generated field of characteristic $0$. Then, $TC^{m} (\overline{X})_{\qit_{\ell}}$ holds $\Leftrightarrow$ the $\ell$-primary torsion subgroup $Br^{m} (\overline{X})[\ell^{\infty}]^{Gal(\overline{k}/k)}$ is finite. 
\end{Prop}
\noindent The finiteness of $Br^{m} (\overline{X})^{Gal(\overline{k}/k)}$ in particular would imply that the cokernel of the restriction map:
\[ Br^{m}(X)[\ell^{\infty}] \to Br^{m}(\overline{X})^{Gal(\overline{k}/k)}[\ell^{\infty}] \]
is at worst finite. Thus, if one believes in the truth of the Tate conjecture, one should expect that the failure of Galois descent for (higher) Brauer classes is at worst finite. An unconditional result in this direction was proved by Colliot-Th\'el\`ene and Skorobogatov for fields of characteristic $0$ when $m=1$; i.e., for the usual Brauer group. Our main result is to extend this Galois descent property to higher Brauer groups:
\begin{Thm}\label{main} Let $X$ be a smooth projective variety over a finitely generated field $k$ and suppose that one of the following holds:
\begin{enumerate}[label=(\alph*)]
\item $k$ has characteristic $0$;
\item $X$ satisfies the standard conjectures (i.e., conjectures $B$, $C$ and $D$ as in \cite{KL})
\end{enumerate}
Then, the cokernel of the map
\[ Br^{m}(X)[\frac{1}{p}] \to Br^{m}(\overline{X})^{Gal(\overline{k}/k)}[\frac{1}{p}] \]
is finite, where $p$ is the exponential characteristic of $k$.
\end{Thm}
\noindent Our plan will be as follows. We first give a proof of Proposition \ref{Tate}, which is fairly routine, given existing results in the literature. The proof of Theorem \ref{main}, on the other hand, will exploit some classical techniques of Deligne from \cite{D1} and \cite{D2}. These allow us to prove some basic degeneracy results for the Hochschild-Serre spectral sequence:
\[ E_{2}^{p,q} = H^{p} (k, H^{q}_{\text{\'et}} (\overline{X}, \qit_{\ell}/\zit_{\ell})) \Rightarrow  H^{p+q}_{\text{\'et}} (X, \qit_{\ell}/\zit_{\ell}).\]
This will involve the introduction of a certain Serre localization which we call the isogeny category, which is the natural context in which to state a degeneracy result for this spectral sequence. Then, we prove a partial (Galois-equivariant) splitting result for the cycle class map with torsion coefficients, and this turns out to be the other important ingredient in the proof of our main result. We do not obtain estimates for the size of the cokernel of 
\[ Br^{m}(X) \to Br^{m}(\overline{X})^{Gal(\overline{k}/k)} \]
as the authors of \cite{CTS} do. This is mostly because any optimal estimate would require that we assume the standard conjectures (even when $k$ has characteristic $0$).
\subsection*{Acknowledgements}
The author would like to thank Jean-Louis Colliot-Th\'el\`ene for pointing out some erroneous statements in an early draft and for suggesting some useful references. The author would also like to thank Bruno Kahn for making some helpful suggestions. Finally, many thanks go to the referee for carefully reading various drafts of this paper and suggesting many important corrections.
\subsection*{Notation}
Throughout this note, we will let $X$ be a smooth projective variety of dimension $d$ over a field $k$, which we will assume to be of characteristic $0$ when necessary. We also let $\overline{k}$ be the separable closure of $k$, $\overline{X}:= X \times_{k} \overline{k}$ and $G_{k}$ the absolute Galois group of $k$. Moreover, for $R = \zit, \qit$ we will use the notation
\[ \hat{R}' := \prod_{\ell \neq \text{char } k} R_{\ell}, \ H^{*}_{\text{\'et}} (-, \hat{R}'(n)) := \prod_{\ell \neq \text{char } k} H^{*}_{\text{\'et}} (-, R_{\ell}(n)),\]
where $\ell$ ranges over primes and where the terms in the above product are the usual $\ell$-adic cohomology groups:
\[ \displaystyle H^{*}_{\text{\'et}} (-, \zit_{\ell}(n)) = \mathop{\lim_{\longleftarrow}}_{m \geq 0} H^{*}_{\text{\'et}} (-, \mu_{\ell^{m}}^{\otimes n}), \ H^{*}_{\text{\'et}} (-, \qit_{\ell}(n)) = H^{*}_{\text{\'et}} (-, \zit_{\ell}(n)) \otimes_{\zit_{\ell}} \qit_{\ell}.\] 
Also, set $\displaystyle \qit/\zit' := \prod_{\ell \neq \text{char } k} \qit_{\ell}/\zit_{\ell}$ and $\displaystyle  H^{*}_{\text{\'et}} (-, \qit/\zit'(n)) := \prod_{\ell \neq \text{char } k} H^{*}_{\text{\'et}} (-, \qit_{\ell}/\zit_{\ell}(n))$, where $\displaystyle H^{*}_{\text{\'et}} (-, \qit_{\ell}/\zit_{\ell}(n)) = \mathop{\lim_{\longrightarrow}}_{m \geq 0} H^{*}_{\text{\'et}} (-, \mu_{\ell^{m}}^{\otimes n}).$
\section{Preliminary Results}

\subsection{The isogeny category}
\begin{Def} Given an additive category $\mathcal{A}$, we define the associated {\em isogeny category of $\A$} to be the Serre localization $\A_{\qit}$ of $\A$ along the isogenies; more precisely, $\A_{\qit}$ is the category whose objects are the same as those of $\A$ and whose morphisms are given by:
\begin{equation} \Hom_{\A_{\qit}} (A, B):=  \Hom_{\A} (A, B) \otimes \qit.\label{iso-hom} \end{equation}
We will say that $\phi \in \Hom_{\A} (A, B)$ is an {\em isogeny} if the corresponding map in $\mathcal{A}_{\qit}$ is an isomorphism.
\end{Def}
\noindent As a matter of convenience, we have the following straightforward lemma:
\begin{Lem}\label{lem-bas} Suppose that $\A$ is an additive category and that $\phi \in \Hom_{\A} (A, B)$ is split-injective (resp., split-surjective) in the isogeny category $\A_{\qit}$; i.e., there exists an integer $M$ and $\psi \in \Hom_{\A} (B, A)$ for which
\[ \psi \circ \phi = M \cdot \text{id}_{A}, \ \text{(resp. } \phi \circ \psi = M \cdot \text{id}_{B}). \]
Then, the kernel (resp., cokernel) of $\phi$ is of finite-exponent.
\end{Lem}
\noindent There is also the following analogue of one of Deligne's decomposition theorems set in the isogeny category:
\begin{Lem}\label{lem-gen} Let $\mathcal{A}$ be an Abelian category and let $D(\A)$ be the derived category of bounded complexes in $\A$. Suppose that there exists $A^{*} \in D(\A)$ and $\phi \in \Hom_{D(\A)} (A^{*}, A^{*}[2])$ for which the map on cohomology
\[ H^{n-i}(\phi^{i}): H^{n-i}(A^{*}) \to H^{n+i}(A^{*}) \]
is an isogeny for all $i >0$. Then, there exists a non-canonical decomposition in $D(\A)_{\qit}$:
\[ A^{*} \cong \bigoplus_{i} H^{i} (A^{*})[-i]. \]
\begin{proof} This is nothing more than \cite{D1} Theorem 1.5 after one realizes that $D(\A_{\qit})$ is the isogeny category of $D(\A)$.
\end{proof}
\end{Lem}
\begin{Cor}\label{main-cor} Retain the notation and assumptions of the previous lemma and let $\Gamma: \A \to \textbf{Ab}$ be a left-exact functor to the category of Abelian groups. Then, for all $i \geq 0$ the edge map
\[ H^{i} (R\Gamma A^{*}) \to \Gamma H^{i} (A^{*})\]
is split-surjective in the isogeny category of Abelian groups; in particular, the cokernel is of finite exponent.
\begin{proof} The first statement follows from Lemma \ref{lem-gen} and \cite{D1} Prop. 1.2. The second statement follows from Lemma \ref{lem-bas}.
\end{proof}
\end{Cor}

\subsection{Some Lefschetz-type Theorems}

Now, let $X$ be a smooth projective variety of dimension $d$ over a field $k$.
\begin{Lem}\label{hard-lef} Suppose that $h \in Pic(X)$ is the class of an ample divisor and $\ell$ is a prime different from the characteristic of $k$.
\begin{enumerate}[label=(\alph*)]
\item\label{l} For all $m\geq 0$, $H^{d-m}_{\text{\'et}} (\overline{X}, \zit_{\ell}) \xrightarrow{\cup h^{m}} H^{d+m}_{\text{\'et}} (\overline{X}, \zit_{\ell}(m))$
is an isogeny in the category of $G_{k}$-modules.
\item\label{no-l} Suppose further that $k$ has characteristic $0$ or that $X$ satisfies the Lefschetz standard conjecture. Then, for all $m\geq 0$ $H^{d-m}_{\text{\'et}} (\overline{X}, \hat{\zit}') \xrightarrow{\cup h^{m}} H^{d+m}_{\text{\'et}} (\overline{X}, \hat{\zit}'(m))$
is an isogeny in the category of $G_{k}$-modules.
\end{enumerate}
\begin{proof} Since $X$ is defined over a finitely generated field, we can assume (by invariance of \'etale cohomology under separably closed extensions, \cite{M} VI Cor. 2.6) that $k$ is some finitely generated field. When $k$ has characteristic $0$, one thus reduces to the case that $k \subset \cit$. In this case, let $H^{*}_{B}(-, \qit)$ denote singular cohomology with $\qit$-coefficients. Then, the classical hard Lefschetz theorem gives an isomorphism:
\[ H^{d-m}_{B} (X_{\cit}, \qit) \xrightarrow{\cup h^{m}} H^{d+m}_{B} (X_{\cit}, \qit)\] 
(ignoring weights). This means that the corresponding map with integral coefficients is an isogeny of Abelian groups. Both statements \ref{l} and \ref{no-l} then follow from the comparison isomorphism between singular and \'etale cohomology. Now suppose that $k$ has positive characteristic. Then, to prove statements \ref{l} and \ref{no-l}, note that by the main result of \cite{D2} we have an isomorphism:
\[H^{d-m}_{\text{\'et}} (\overline{X}, \qit_{\ell}) \xrightarrow{\cup h^{m}} H^{d+m}_{\text{\'et}} (\overline{X}, \qit_{\ell}(m)). \]
Since $H^{*}_{\text{\'et}} (\overline{X}, \zit_{\ell})$ is a finitely generated $\zit_{\ell}$ module, it follows that the corresponding map with $\zit_{\ell}$ coefficients is an isogeny of $G_{k}$-modules. To obtain the corresponding statement for $H^{*}_{\text{\'et}} (\overline{X}, \hat{\zit}')$, we use the Lefschetz standard conjecture to obtain a correspondence
\[ \Lambda_{m} \in CH^{d-m} (X \times X) \otimes \qit\]
for which $(\cup h^{m})\circ \Lambda_{m} = \text{id}_{H^{d+m}_{\text{\'et}} (\overline{X}, \qit_{\ell}(m))}$ and $\Lambda_{m}\circ(\cup h^{m}) = \text{id}_{H^{d-m}_{\text{\'et}} (\overline{X}, \qit_{\ell})}$ for all primes $\ell \neq \text{char } k$. Of course, we can find some integer $M$ (not depending on $\ell$) for which $\tilde{\Lambda}_{m}:= M\cdot \Lambda_{m} \in CH^{d-m} (X \times X)$ and since 
\[ \tilde{\Lambda}_{m}(H^{d+m}_{\text{\'et}} (\overline{X}, \zit_{\ell}(m))) \subset H^{d-m}_{\text{\'et}} (\overline{X}, \zit_{\ell}) \]
for all primes $\ell \neq \text{char } k$, it follows that $(\cup h^{m})\circ \tilde{\Lambda}_{m} = M\cdot\text{id}_{H^{d+m}_{\text{\'et}} (\overline{X}, \zit_{\ell}(m))}$ and $\tilde{\Lambda}_{m}\circ(\cup h^{m}) = M\cdot\text{id}_{H^{d-m}_{\text{\'et}} (\overline{X}, \zit_{\ell})}$.
Since $M$ is invertible in $\zit_{\ell}$ for all but finitely many primes $\ell$, this means that for all such primes $\cup h^{m}: H^{d-m}_{\text{\'et}} (\overline{X}, \zit_{\ell}) \to H^{d+m}_{\text{\'et}} (\overline{X}, \zit_{\ell}(m))$ is an isomorphism. Taking the product over all $\ell$, we deduce that
\[\cup h^{m}: H^{d-m}_{\text{\'et}} (\overline{X}, \hat{\zit}') \to H^{d+m}_{\text{\'et}} (\overline{X}, \hat{\zit}'(m))\]
 is an isogeny of $G_{k}$-modules, which gives statement \ref{no-l} in positive characteristic.
\end{proof}
\end{Lem}
\begin{Rem} The reader should note that the role of the Lefschetz standard conjecture in characteristic $p>0$ is to ensure that the degree of the isogeny 
\[\cup h^{m}: H^{d-m}_{\text{\'et}} (\overline{X}, \zit_{\ell}) \to H^{d+m}_{\text{\'et}} (\overline{X}, \zit_{\ell}(m))\]
does not depend on $\ell$. The author is not sure how to prove this in the absence of the Lefschetz standard conjecture.
\end{Rem}
\begin{Cor}\label{hard-lef-cor} Retain the set-up of Lemma \ref{hard-lef}.
\begin{enumerate}[label=(\alph*)]
\item\label{l-2} For all $m\geq 0$, $H^{d-m}_{\text{\'et}} (\overline{X}, \qit_{\ell}/\zit_{\ell}) \xrightarrow{\cup h^{m}} H^{d+m}_{\text{\'et}} (\overline{X}, \qit_{\ell}/\zit_{\ell}(m))$
is an isogeny in the category of $G_{k}$-modules.
\item\label{no-l-2} Suppose further that $k$ has characteristic $0$ or that $X$ satisfies the Lefschetz standard conjecture. Then, for all $m\geq 0$ $H^{d-m}_{\text{\'et}} (\overline{X}, \qit/\zit') \xrightarrow{\cup h^{m}} H^{d+m}_{\text{\'et}} (\overline{X}, \qit/\zit(m)')$
is an isogeny in the category of $G_{k}$-modules.
\end{enumerate}
\begin{proof}
For statement \ref{l-2}, there is a commutative diagram with rows exact:
\[\begin{tikzcd}  0 \arrow{r} & H^{d-m}_{\text{\'et}} (\overline{X}, \zit_{\ell}) \otimes \qit_{\ell}/\zit_{\ell} \arrow{r} \arrow{d}{\cup h^{d-m}} & H^{d-m}_{\text{\'et}} (\overline{X}, \qit_{\ell}/\zit_{\ell}) \arrow{r} \arrow{d}{\cup h^{d-m}} & H^{d-m+1}_{\text{\'et}} (\overline{X}, \zit_{\ell})[\ell^{\infty}] \arrow{d}{\cup h^{d-m}} \arrow{r} & 0\\ 
0 \arrow{r} & H^{d+m}_{\text{\'et}} (\overline{X}, \zit_{\ell}) \otimes \qit_{\ell}/\zit_{\ell} \arrow{r}  & H^{d+m}_{\text{\'et}} (\overline{X}, \qit_{\ell}/\zit_{\ell}) \arrow{r} & H^{d+m+1}_{\text{\'et}} (\overline{X}, \zit_{\ell})[\ell^{\infty}] \arrow{r} & 0
\end{tikzcd}\]
(suppressing weights). Since $H^{*}_{\text{\'et}} (\overline{X}, \zit_{\ell})$ is a finitely generated $\zit_{\ell}$-module, both the rightmost terms are finite and, hence, vanish in the isogeny category. So, to prove that the middle vertical arrow is an isogeny, it suffices to prove that the left vertical arrow is an isogeny. But this latter follows directly from Lemma \ref{hard-lef} \ref{l} by tensoring with $\qit_{\ell}/\zit_{\ell}$. To prove the statement of \ref{no-l-2}, note that there is an identical diagram with coefficients in $\qit/\zit'$ and with the rightmost terms $H^{d-m+1}_{\text{\'et}} (\overline{X}, \hat{\zit}')_{tors}$ and $H^{d+m+1}_{\text{\'et}} (\overline{X}, \hat{\zit}')_{tors}$. By the comparison isomorphism with singular cohomology, $H^{2m+1}_{\text{\'et}} (\overline{X}, \hat{\zit}'(m))_{tors}$ is finite (since singular cohomology with integral coefficients is finitely generated), when $k$ is of characteristic $0$. On the other hand, when $\text{char } k >0$ the finiteness of this group follows from the main result of \cite{G}. Since finite torsion groups vanish in the isogeny category of Abelian groups (and, hence, of $G_{k}$-modules), the argument from statement \ref{l-2} works mutatis mutandis (using Lemma \ref{hard-lef} \ref{no-l} this time).
\end{proof}
\end{Cor}
\begin{Cor}\label{natural} Retain the assumptions of Lemma \ref{hard-lef}.
\begin{enumerate}[label=(\alph*)]
\item\label{l-3} For all $n$ and $m$, the cokernel of the natural map
\begin{equation} H^{m}_{\text{\'et}} (X, \qit_{\ell}/\zit_{\ell}(n)) \to H^{m}_{\text{\'et}} (\overline{X}, \qit_{\ell}/\zit_{\ell}(n))^{G_{k}}\label{weak} \end{equation}
is finite.
\item\label{no-l-3} Suppose further that $k$ has characteristic $0$ or that $X$ satisfies the Lefschetz standard conjecture. Then, for all $n$ and $m$, the cokernel of the natural map
\begin{equation} H^{m}_{\text{\'et}} (X, \qit/\zit'(n)) \to H^{m}_{\text{\'et}} (\overline{X}, \qit/\zit'(n))^{G_{k}}\label{weak-2} \end{equation}
is finite.
\end{enumerate}

\begin{proof} For both statements, we first observe that the cokernel is of finite exponent, which follows for (\ref{weak}) (resp., for (\ref{weak-2})) by Corollary \ref{hard-lef-cor} \ref{l-2} and Corollary \ref{main-cor} (resp., by Corollary \ref{hard-lef-cor} \ref{no-l-2} and Corollary \ref{main-cor}). The desired statement then follows from the lemma below, which is a straightforward group-theoretic fact (cf, \cite{CTS} \S 1.2).
\end{proof} 
\end{Cor}
\begin{Lem}\label{nat-lem} Suppose that $A$ is an Abelian group of cofinite type; i.e., of the form $(\qit_{\ell}/\zit_{\ell})^{r} \oplus F$ (or $(\qit/\zit')^{r} \oplus F$), where $r \geq 0$ and $F$ is finite. Then, any finite-exponent subquotient of $A$ is finite.
\end{Lem}

\begin{Rem} We remark that Corollary \ref{natural} is only interesting for $m = 2n$ since for all others, the right hand groups are already finite by the results of \cite{CTR}. 
\end{Rem}
\section{Higher Brauer groups}
\subsection{\'Etale motivic cohomology}
\noindent For $n \geq 0$ an integer, let $\tau$ denote either the Zariski or the \'etale topology. To define the corresponding motivic cohomology groups, let $z^{n} (-, *)_{\tau}$ be the complex of sheaves for the $\tau$ topology on $X$ given by Bloch’s cycle complex \cite{Bl}. Then, for any Abelian group $A$ consider the {\em cycle complex}:
\[ A_{X}(n)_{\tau} := (z^{n} (-, *)_{\tau} \otimes A)[-2n]. \]
Then, the motivic and \'etale motivic cohomology groups are given by the corresponding hypercohomology groups:
\[ \begin{split}
H^{p}_{M} (X, A(n)) & = \mathbb{H}^{p}_{\text{Zar}} (X, A_{X}(n)_{\text{Zar}})\\
H^{p}_{L} (X, A(n)) & = \mathbb{H}^{p}_{\text{\'et}} (X, A_{X}(n)_{\text{\'et}}).
\end{split}
\]
The motivic cohomology groups are nothing more than Bloch's higher Chow groups, $H^{p}_{M} (X, \zit(n)) = CH^{n} (X, 2n-p)$, since Bloch's cycle complex satisfies Zariski descent (see \cite{Bl} \S 3). As per convention, we denote the Lichtenbaum Chow group by $CH^{n}_{L} (X) =  H^{2n}_{L} (X, \zit(n))$. 
There are also the following important facts due to numerous authors:
\begin{Thm}\label{props} With the above notation,
\begin{enumerate}[label=(\alph*)]
\item\label{q-coeff} $H^{p}_{M} (X, \qit(n)) \cong H^{p}_{L} (X, \qit(n))$;
\item\label{iso-1} the cycle class map $c_{m}^{p, n}: H^{p}_{L} (X, \zit/m(n)) \to  H^{p}_{\text{\'et}} (X, \mu_{m}^{\otimes n})$ is an isomorphism for $m$ invertible in $k$;
\item\label{short-exact-mot} there exists a natural short exact sequence:
\[ 0 \to H^{p}_{L} (X, \zit(n)) \otimes \qit/\zit' \to H^{p}_{\text{\'et}} (X, \qit/\zit'(n)) \to H^{p+1}_{L} (X, \zit(n))_{tors}[\frac{1}{p}] \to 0, \]
where the first non-zero arrow is the cycle class map and $p$ is the exponential characteristic of $k$;
\item\label{tors-mot} $H^{p}_{L} (X, \zit(n))$ is torsion for $p>2n$.
\end{enumerate}
\begin{proof} The statement \ref{q-coeff} is part of \cite{K} Th\'eor\`eme 2.6. The statement of \ref{iso-1} is the main result of Geisser and Levine's paper \cite{GL2}. The short exact sequence \ref{short-exact-mot} comes from the coefficient theorem. Statement \ref{tors-mot} follows from \ref{q-coeff} and the fact that the higher Chow groups $H^{p}_{M} (X, \zit(m)) \cong CH^{m} (X, 2m-p)$ vanish for $p>2m$. 
\end{proof}
\end{Thm}
\noindent Now, we use the following notation for the image of the integral cycle class map:
\begin{Not} Let $A^{n} (\overline{X})$ denote the image of the total cycle class map:
\[c^{n, n}:= \mathop{\lim_{\longleftarrow}}_{(m, \text{char } k) = 1} c_{m}^{n, n}:  CH^{n}_{L} (\overline{X}) \to H^{2n}_{\text{\'et}} (\overline{X}, \hat{\zit}'(n)).\]
\end{Not}
\noindent Note that by Theorem \ref{props} \ref{q-coeff} $A^{n} (\overline{X}) \otimes \qit$ coincides with the image of the usual cycle class map:
\[ CH^{n} (\overline{X}) \otimes \qit \to H^{2n}_{\text{\'et}} (\overline{X}, \hat{\qit}'(n)).  \]
Moreover, we denote by $Br^{n} (X) := H^{2n+1}_{L} (X, \zit(n))$ the {\em $n^{th}$ higher Brauer group}, which is a torsion group by \ref{tors-mot} of this theorem. Finally, there is the following consequence of the above theorem:
\begin{Cor}\label{cor-desired-short} Let $k$ be a field of exponential characteristic $p$ and let $X$ be a smooth projective variety over $k$. Then, there is a natural short exact sequence of $G_{k}$-modules:
\begin{equation} 0 \to A^{n} (\overline{X}) \otimes \qit/\zit' \to H^{2n}_{\text{\'et}} (\overline{X}, \qit/\zit'(n)) \to  Br^{n}(\overline{X})[\frac{1}{p}] \to 0,\label{desired-short}\end{equation}
where the left non-zero arrow is the cycle class map.
\begin{proof} This essentially follows from Theorem \ref{props} \ref{short-exact-mot}. Indeed, this latter gives a short exact sequence:
\begin{equation} 0 \to CH^{n}_{L} (\overline{X}) \otimes \qit/\zit' \to H^{2n}_{\text{\'et}} (\overline{X}, \qit/\zit'(n)) \to Br^{n} (\overline{X})[\frac{1}{p}] \to 0.\label{short} \end{equation}
To see that one can replace $CH^{n}_{L} (\overline{X})$ with $A^{n} (\overline{X})$, note the obvious factorization:
\[\begin{tikzcd}
CH^{n}_{L} (\overline{X}) \otimes \qit/\zit' \arrow[two heads]{d} \arrow[hook]{r} & H^{2n}_{\text{\'et}} (\overline{X}, \qit/\zit'(n)) \\
A^{n} (\overline{X}) \otimes \qit/\zit' \arrow{r} & H^{2n}_{\text{\'et}} (\overline{X}, \hat{\zit}'(n)) \otimes_{\hat{\zit}'} \qit/\zit'. \arrow[hook]{u}
\end{tikzcd} \]
Here, the top horizontal arrow is the left non-zero arrow from (\ref{short}). Since this latter is injective, so is the left vertical arrow. Hence, the left vertical arrow is an isomorphism, as desired.
\end{proof}
\end{Cor}

\subsection{Absolute Hodge classes}
\noindent In this subsection, we will assume that $k$ is a finitely generated field of characteristic $0$. For $R = \zit, \qit$ we suppress the $'$ notation from earlier and write $\hat{R} = \hat{R}'$ and $H^{2m}_{\text{\'et}} (\overline{X}, \hat{R}(m)) = H^{2m}_{\text{\'et}} (\overline{X}, \hat{R}'(m))$. We also write $H^{*}_{B}(-, R(n))$ for singular cohomology with coefficients in $R(n)$ for $R= \zit, \qit$ (note that, as usual, $R(n) = R(1)^{\otimes n}$ with $R(1) = R(2\pi i) \subset \cit$).\\
\indent Now, we would like to recall the notion of absolute Hodge classes, first introduced by Deligne in \cite{D3}. Indeed, for any embedding of fields $\sigma : k \hookrightarrow \cit$, we consider the corresponding comparison isomorphism of \'etale cohomology with singular cohomology:
\[ \text{comp}_{\sigma}^{m,n}: H^{m}_{\text{\'et}} (\overline{X}, \hat{\qit}(n)) \xrightarrow{\cong} H^{m}_{\text{\'et}} (X \times_{\sigma} \cit, \hat{\qit}(n)) \to H^{m}_{B} (X \times_{\sigma} \cit, \qit(n)) \otimes_{\qit}  \hat{\qit},\]
where the first isomorphism comes from the invariance of \'etale cohomology under algebraically closed extensions. Moreover, for each such embedding, we consider the corresponding subspace of Hodge classes:
\[ H^{n, n}(X \times_{\sigma} \cit, \qit) \subset H^{2n}_{B} (X \times_{\sigma} \cit, \qit(n)) \subset   H^{2n}_{B} (X \times_{\sigma} \cit, \qit(n)) \otimes_{\qit}  \hat{\qit}.\]
Then, we recall that the subspace of absolute Hodge classes is defined to be:
\[  AH^{n} (\overline{X}) := \bigcap_{\sigma: k \hookrightarrow \cit} (\text{comp}_{\sigma}^{2n,n})^{-1} (H^{n, n}(X \times_{\sigma} \cit, \qit)).\]

\noindent There are the following important observations about absolute Hodge cycles that we will need in the sequel:
\begin{Lem}\label{Hodge} Suppose that $X$ is a smooth projective variety over $k$. Then, for all non-negative integers $n$, we have:
\begin{enumerate}[label=(\alph*)]
\item\label{Hodge-cont} $A^{n} (\overline{X}) \otimes \qit \subset AH^{n} (\overline{X})$;
\item\label{Hodge-fin} $AH^{n} (\overline{X})$ is a $G_{k}$-stable subspace of $H^{2n}_{\text{\'et}} (\overline{X}, \hat{\qit}(n))$;
\item\label{Hodge-non} There exists a non-degenerate $G_{k}$-invariant pairing:
\[ (,): H^{2n}_{\text{\'et}} (\overline{X}, \hat{\qit}(n)) \otimes H^{2n}_{\text{\'et}} (\overline{X}, \hat{\qit}(n)) \to \hat{\qit}\]
for which the restriction of $(,)$ to $AH^{n} (\overline{X})$ is $\qit$-valued and non-degenerate.
\end{enumerate}
\begin{proof} Statement \ref{Hodge-cont} follows from the fact that the comparison isomorphism above is compatible with cycle class maps and that the image of the cycle class map for singular cohomology lies in the subspace of Hodge classes. Statement \ref{Hodge-fin} is easy to see from the definition. Finally, statement \ref{Hodge-non} follows from Proposition 22 of \cite{CS} (see also \cite{A} Prop. 3.3).
\end{proof}

\end{Lem}

\subsection{Some splitting results}
\noindent For the applications to follow, we will need some ad hoc results on the split-injectivity (as a map of Galois modules) of the cycle class map with torsion coefficients. 
\begin{Lem}\label{above} Suppose that $k$ is a finitely generated field of characteristic $0$ and that $X$ is a smooth projective variety over $k$. 
\begin{enumerate}[label=(\alph*)]
\item\label{above-1} $AH^{n} (\overline{X})$ is a finite rank $\qit$-vector space. In particular, there exists a finite extension $k'/k$ for which the action of $G_{k'}$ on $AH^{n} (\overline{X})$ is trivial.
\item\label{above-2} $A^{n} (\overline{X})/tors$ is a finitely generated subgroup of the $\qit$-vector space $AH^{n} (\overline{X})$.
\end{enumerate}
\begin{proof} The first statement of \ref{above-1} is clear by definition of $AH^{n} (\overline{X})$. The second statement of \ref{above-1} then follows from the first and the fact that a pro-finite group acting continuously on a finite rank $\qit$-vector space acts via a finite quotient. For \ref{above-2}, we observe that by Lemma \ref{Hodge} \ref{Hodge-cont}, it follows that $A^{n} (\overline{X})/tors \subset  AH^{n} (\overline{X})$. Now, we note that from item \ref{above-1} of this lemma, $A^{n} (\overline{X})/tors$ is of finite rank, but this is insufficient to show that it is finitely generated. To this last end, let 
\[ H^{n,n} (X \times_{\sigma} \cit, \zit) \subset H^{2n}_{B} (X \times_{\sigma} \cit, \zit(n))/tors \]
be the subgroup of integral Hodge classes of degree $n$ on $X \times_{\sigma} \cit$. Then, we observe that $A^{n} (\overline{X})/tors$ lies in the subgroup
\[ \bigcap_{\sigma: k \hookrightarrow \cit} (\text{comp}_{\sigma}^{2n,n})^{-1} (H^{n, n}(X \times_{\sigma} \cit, \zit)) \subset AH^{n} (\overline{X})\]
This subgroup is a lattice in $AH^{n} (\overline{X})$, which means that $A^{n} (\overline{X})/tors$ is finitely generated, as desired.
\end{proof}
\end{Lem}

\begin{Cor}\label{cor-above} With the set-up of Lemma \ref{above}, the cokernel of the inclusion
\begin{equation} A^{n} (\overline{X})^{G_{k}} \otimes \qit/\zit \hookrightarrow (A^{n} (\overline{X})\otimes \qit/\zit)^{G_{k}}\label{arrow} \end{equation}
is a group of finite exponent.
\begin{proof} There is the obvious short-exact sequence of $G_{k}$-modules:
\[ 0 \to A^{n} (\overline{X})/\text{tors} \to A^{n} (\overline{X})\otimes \qit \to A^{n} (\overline{X})\otimes \qit/\zit \to 0. \] 
Applying $G_{k}$-invariants gives an exact sequence:
\[ 0 \to  A^{n} (\overline{X})^{G_{k}} \otimes \qit/\zit \to (A^{n} (\overline{X})\otimes \qit/\zit)^{G_{k}} \to H^{1} (G_{k}, A^{n} (\overline{X})/\text{tors}). \]
So, it suffices to show that $H^{1} (G_{k}, A^{n} (\overline{X})/\text{tors})$ is a group of finite exponent. Using the inflation-restriction exact sequence in group cohomology, we are free to pass to a finite extension of $k$. By Lemma \ref{above}, we can assume that $G_{k}$ acts trivially on $A^{n} (\overline{X})/\text{tors}$. Since $A^{n} (\overline{X})/\text{tors}$ is also torsion-free, we have $H^{1} (G_{k}, A^{n} (\overline{X})/\text{tors}) = 0$, which gives the desired result.
\end{proof}
\end{Cor}

\noindent The result below is an elementary lemma in representation theory and is included for convenience.

\begin{Lem}\label{split} Suppose that $W$ is a finite rank free module over a complete discrete valuation ring $R$, that $G$ is a group acting $R$-linearly on $W$ and that there is a non-degenerate $G$-invariant pairing
\[ (,): W \otimes_{R} W \to R \]
and a $G$-stable submodule $V \subset W$ for which $(,)$ restricted to $V$ is also non-degenerate. Then, the inclusion map $i:V \hookrightarrow W$ is split-injective in the isogeny category of $G$-modules. 
\begin{proof} The pairing above gives a map of $G$-modules, $\phi: W \to W^{\vee} := \Hom_{R} (W, R)$. Moreover, the non-degeneracy assumption of $\restr{(,)}{V}$ ensures that the composition of maps of $G$-modules:
\[ \psi: V \xhookrightarrow{i} W \xrightarrow{\phi} W^{\vee} \xrightarrow{res}  V^{\vee}\]
is injective. Indeed, $\psi$ is nothing more than the map $\alpha \mapsto (\alpha, -)$. Thus, for rank reasons, $\psi$ is an isogeny of $G$-modules, and a left-inverse for $i$ is then given by the composition
\[ W \xrightarrow{\phi} W^{\vee} \xrightarrow{res}  V^{\vee} \xrightarrow[\cong]{\psi^{-1}} V \]
where $V^{\vee} \xrightarrow[\cong]{\psi^{-1}} V$ is the inverse isogeny. So, $i: V \to W$ is split-injective in the isogeny category of $G$-modules, as desired.
\end{proof}
\end{Lem}

\begin{Rem}\label{rem-conv} For the applications below, the main interest is the case that $R = \zit_{\ell}$ for some prime $\ell$, $G$ is the absolute Galois group of a field $k$ of characteristic $0$, $W = H^{2n} (\overline{X}, \zit_{\ell}(n))/tors$ (for some smooth projective variety $X$ over $k$), $V = M \otimes \zit_{\ell}$ with $M$ a $G$-stable finitely generated subgroup in $AH^{n} (\overline{X})$ and $(,)$ is the pairing in Lemma \ref{Hodge} \ref{Hodge-non}. Since $(\alpha, \beta) \subset \qit$ for $\alpha, \beta \in AH^{n} (\overline{X})$, it is easy to see that the degree of the isogeny $\psi: V \to V^{\vee}$ in the proof of Lemma \ref{split} is the discriminant of the pairing $(,)$ on $M$. Thus, for all but finitely many $\ell$, $\psi$ is an isomorphism and the inclusion $V \hookrightarrow W$ is split-injective (without passing to the isogeny category). 
\end{Rem}

\begin{Lem}\label{split-2} Suppose that $X$ is a smooth projective variety over $k$ and that $k$ is a field of characteristic $0$ or that $X$ satisfies the standard conjectures (i.e., conjectures $B$, $C$ and $D$ as in \cite{KL}). Then, for all non-negative integers $n$, the natural inclusion
\begin{equation} A^{n} (\overline{X})/tors \otimes_{\zit} \hat{\zit}' \hookrightarrow H^{2n}_{\text{\'et}} (\overline{X}, \hat{\zit}'(n))/tors\label{goal-inj} \end{equation}
is split-injective in the isogeny category of $G_{k}$-modules.
\begin{proof} 
First we assume that $k$ has characteristic $0$. By Lemma \ref{above}, we can assume without loss of generality that $G_{k}$ acts trivially on $AH^{n} (\overline{X})$. Indeed, let $k'/k$ be a finite Galois extension for which $G_{k'}$ acts trivially on $AH^{n} (\overline{X})$ and observe that if a splitting $\rho$ exists which is $G_{k'}$-invariant, then one can set
\[ \rho' := \sum_{h \in H} h^{*}\rho, \]
where $H = G_{k}/G_{k'}$ to obtain a $G_{k}$-invariant splitting.\\ 
\indent Assuming that $G_{k}$ acts trivially on $AH^{n} (\overline{X})$, let $M = A^{n} (\overline{X})/tors$ (by Lemma \ref{above}, this is a finitely generated subgroup of $AH^{n} (\overline{X})$). Also, let $\ell$ be a prime and set $V_{\ell} := M \otimes \zit_{\ell}$. From Remark \ref{rem-conv}, the inclusion $V_{\ell} \hookrightarrow H^{2n} (\overline{X}, \zit_{\ell}(n))$ is split-injective for all but finitely many $\ell$ and in general split-injective in the isogeny category of $G_{k}$-modules. Thus, taking a product over $\ell$, we deduce that the inclusion
\begin{equation} M \otimes \hat{\zit} \hookrightarrow H^{2n}_{\text{\'et}} (\overline{X}, \hat{\zit}(n))\label{spl} \end{equation}
is also split-injective in the isogeny category of $G_{k}$-modules. This gives the result in the case of characteristic $0$.\\
\indent On the other hand, suppose that $X$ satisfies the standard conjectures. By standard conjecture $D$, for any ample divisor $h \in Pic(\overline{X})$ the non-degenerate $G_{k}$-invariant pairing
\[ (,): H^{2n}_{\text{\'et}} (\overline{X}, \hat{\qit}'(n)) \otimes H^{2n}_{\text{\'et}} (\overline{X}, \hat{\qit}'(n)) \to \hat{\qit}', \ (\alpha, \beta) = \alpha\cup \beta \cup h^{d-2n}\]
for $2n \leq d$ induces a non-degenerate $\qit$-valued pairing on $A^{n} (\overline{X}) \otimes \qit$. Additionally, by the Lefschetz standard conjecture, the inverse $\Lambda_{n}$ to the isomorphism:
\begin{equation} H^{2n}_{\text{\'et}} (\overline{X}, \hat{\qit}'(n)) \xrightarrow{\cup h^{d-2n}} H^{2(d-n)}_{\text{\'et}} (\overline{X}, \hat{\qit}'(d-n))\label{cup-h}\end{equation}
is algebraic (for $2n \leq d$); note that $\Lambda_{n}$ is $G_{k}$-equivariant since (\ref{cup-h}) is. Thus, the corresponding $G_{k}$-invariant pairing for $2n>d$:
\[ H^{2n}_{\text{\'et}} (\overline{X}, \hat{\qit}'(n)) \otimes H^{2n}_{\text{\'et}} (\overline{X}, \hat{\qit}'(n)) \to \hat{\qit}', \ (\alpha, \beta) = \Lambda_{n}(\alpha)\cup \beta \]
is non-degenerate and induces a non-degenerate $\qit$-valued pairing on $A^{n} (\overline{X}) \otimes \qit$. The same rationale as in Remark \ref{rem-conv} (with $M = A^{n} (\overline{X})$) then gives the desired splitting result.
\end{proof}
\end{Lem}

\begin{Cor}\label{splitting} With the set-up of Lemma \ref{split-2}, the natural inclusion
\[ A^{n} (\overline{X}) \otimes \qit/\zit' \hookrightarrow H^{2n}_{\text{\'et}} (\overline{X}, \qit/\zit'(n))   \]
is split-injective in the isogeny category of $G_{k}$-modules.
\begin{proof} The above inclusion factors as 
\[ A^{n} (\overline{X}) \otimes \qit/\zit' \hookrightarrow H^{2m}_{\text{\'et}} (\overline{X}, \hat{\zit}'(m)) \otimes \qit/\zit' \hookrightarrow H^{2n}_{\text{\'et}} (\overline{X}, \qit/\zit'(n)). \]
By Lemma \ref{split-2}, the first inclusion is split-injective in the isogeny category of $G_{k}$-modules. On the other hand, there is the short exact sequence of $G_{k}$-modules:
\[ 0 \to H^{2m}_{\text{\'et}} (\overline{X}, \hat{\zit}'(m)) \otimes \qit/\zit' \to H^{2m}_{\text{\'et}} (\overline{X}, \qit/\zit'(m)) \to H^{2m+1}_{\text{\'et}} (\overline{X}, \hat{\zit}'(m))_{tors} \to 0. \]
By the comparison isomorphism with singular cohomology, $H^{2m+1}_{\text{\'et}} (\overline{X}, \hat{\zit}'(m))_{tors}$ is finite (since singular cohomology with integral coefficients is finitely generated), when $k$ is of characteristic $0$. On the other hand, when $\text{char } k >0$ the finiteness of this group follows from the main result of \cite{G}. Since finite torsion groups vanish in the isogeny category of Abelian groups (and, hence, of $G_{k}$-modules), it follows that $H^{2m}_{\text{\'et}} (\overline{X}, \hat{\zit}'(m)) \otimes \qit/\zit' \to H^{2m}_{\text{\'et}} (\overline{X}, \qit/\zit'(m))$ is an isomorphism in the isogeny category of $G_{k}$-modules. This gives the desired result.
\end{proof}
\end{Cor}

\section{Proofs of main results}
\subsection{Proof of Proposition \ref{Tate}}
\noindent (Proof of $\Leftarrow$) Consider the short exact sequence of $G_{k}$-modules:
\begin{equation} 0 \to A^{m} (\overline{X}) \otimes \qit_{\ell}/\zit_{\ell} \to H^{2m}_{\text{\'et}} (\overline{X}, \qit_{\ell}/\zit_{\ell}(m)) \to Br^{m} (\overline{X})[\ell^{\infty}] \to 0.\label{exact-1} \end{equation}
Applying the $G_{k}$-invariants functor to (\ref{exact-1}), we obtain an exact sequence:
\[ 0 \to (A^{m} (\overline{X}) \otimes \qit_{\ell}/\zit_{\ell})^{G_{k}} \xrightarrow{\phi} H^{2m}_{\text{\'et}} (\overline{X}, \qit_{\ell}/\zit_{\ell}(m))^{G_{k}} \to Br^{m} (\overline{X})[\ell^{\infty}]^{G_{k}}.\]
Suppose that $Br^{m} (\overline{X})[\ell^{\infty}]^{G_{k}}$ is finite, so that $\text{coker } \phi$ is also. Then, there is a commutative diagram with rows exact:
\[ \begin{tikzcd}
& A^{m} (\overline{X})^{G_{k}} \otimes \zit_{\ell} \arrow{r} \arrow{d} & A^{m} (\overline{X})^{G_{k}} \otimes \qit_{\ell} \arrow{r} \arrow{d} & A^{m} (\overline{X})^{G_{k}} \otimes \qit_{\ell}/\zit_{\ell}\arrow{d}{\phi'} \arrow{r} & 0\\
0 \arrow{r} & H^{2m}_{\text{\'et}} (\overline{X}, \zit_{\ell}(m))^{G_{k}}/tors \arrow{r}  & H^{2m}_{\text{\'et}} (\overline{X}, \qit_{\ell}(m))^{G_{k}} \arrow{r} & H^{2m}_{\text{\'et}} (\overline{X}, \qit_{\ell}/\zit_{\ell}(m))^{G_{k}}.
\end{tikzcd} \]
Now, let $C:= \text{coker } \{ A^{m} (\overline{X})^{G_{k}} \otimes \zit_{\ell} \to H^{2m}_{\text{\'et}} (\overline{X}, \zit_{\ell}(m))^{G_{k}}/tors \}$. Then, a diagram chase shows that (since $\phi'$ is injective) there is an injective map
\begin{equation} C  \otimes_{\zit_{\ell}} \qit_{\ell}/\zit_{\ell} \hookrightarrow \text{coker } \phi'. \label{hook}\end{equation} 
We first would like to show that $C \otimes_{\zit_{\ell}} \qit_{\ell}/\zit_{\ell}$ is a group of finite exponent. To this end, (\ref{hook}) implies that it suffices to show $\text{coker } \phi'$ is a group of finite exponent. Now, we have the composition
\[ \phi': A^{m} (\overline{X})^{G_{k}} \otimes \qit_{\ell}/\zit_{\ell} \xhookrightarrow{(\ref{arrow})} (A^{m} (\overline{X}) \otimes \qit_{\ell}/\zit_{\ell})^{G_{k}} \xhookrightarrow{\phi} H^{2m}_{\text{\'et}} (\overline{X}, \qit_{\ell}/\zit_{\ell}(m))^{G_{k}}.\]
Note that the cokernel of (\ref{arrow}) is of finite exponent by Corollary \ref{cor-above}. Since $\text{coker } \phi$ was already noted to be finite, it follows that $\text{coker } \phi'$ is of finite exponent, as desired. Now, since $C$ is a finitely generated $\zit_{\ell}$-module, it follows that $C \otimes_{\zit_{\ell}} \qit_{\ell}/\zit_{\ell} = 0$ and hence that $C$ is finite. Thus, $C \otimes_{\zit_{\ell}} \qit_{\ell} = 0$, which implies that $TC^{m} (\overline{X})_{\qit_{\ell}}$ holds.\\
(Proof of $\Rightarrow$) Suppose conversely that the Tate conjecture $TC^{m} (\overline{X})_{\qit_{\ell}}$ holds. Then, the proof proceeds as in \cite{SZ} Prop. 2.5. Indeed, let
\[ T_{\ell}(Br^{m} (\overline{X})) :=  \Hom(\qit_{\ell}/\zit_{\ell}, Br^{m} (\overline{X})) =\mathop{\lim_{\longleftarrow}}_{r\geq 0} Br^{m}(\overline{X})[\ell^{r}].\]
This group is torsion-free by \cite{RS} Theorem 1.3. Then, by Corollary \ref{splitting} there is a splitting of $G_{k}$-modules 
\begin{equation} H^{2m}_{\text{\'et}} (\overline{X},\qit_{\ell}(m)) \cong A^{m}(\overline{X}) \otimes \qit_{\ell} \oplus T_{\ell}(Br^{m} (\overline{X})) \otimes_{\zit_{\ell}} \qit_{\ell}.\label{split-exact} \end{equation} 
Moreover, since $A^{m} (\overline{X}) \otimes \qit$ is a finite rank $\qit$-vector space by Lemma \ref{above}, it follows from the Tate conjecture there exists some finite Galois extension $k'/k$ for which 
\[  A^{m}(\overline{X})^{G_{k'}} \otimes \qit_{\ell} \xrightarrow{\cong} H^{2m}_{\text{\'et}} (\overline{X},\qit_{\ell}(m))^{G_{k'}}  \] 
This implies that $T_{\ell}(Br^{m} (\overline{X}))^{G_{k'}}\otimes_{\zit_{\ell}} \qit_{\ell} = 0$, which means that $T_{\ell}(Br^{m} (\overline{X})^{G_{k'}}) = T_{\ell}(Br^{m} (\overline{X}))^{G_{k'}} = 0$ (since $T_{\ell}(Br^{m} (\overline{X}))$ is torsion-free). 
Since $Br^{m} (\overline{X})^{G_{k'}}[\ell^{\infty}]$ is of cofinite type, it now follows easily that $Br^{m} (\overline{X})^{G_{k'}}[\ell^{\infty}]$ must be finite.

\begin{Rem} The main difficulty in extending this proof to characteristic $p>0$ is that there is no notion of absolute Hodge classes in this case and so even with the Tate conjecture for $X$, it is not clear that the splitting of $G_{k}$-modules (\ref{split-exact}) holds. If one instead assumes the Tate conjecture for all products $X^{n}$, then the results of \cite{Mo} (for instance) show that $H^{2m}_{\text{\'et}} (\overline{X},\qit_{\ell}(m))$ is semi-simple as a $G_{k}$-module. This would then imply (\ref{split-exact}). 
\end{Rem}

\subsection{Proof of Theorem \ref{main}}

\noindent Note that there is a commutative diagram with rows exact:
\[\begin{tikzcd}
& CH^{m}_{L} (X) \otimes \qit/\zit' \arrow{r} \arrow{d} & H^{2m}_{\text{\'et}} (X, \qit/\zit'(m)) \arrow{r} \arrow{d} & Br^{m} (X)[\frac{1}{p}] \arrow{r} \arrow{d} &0\\
0 \arrow{r} & (A^{m} (\overline{X}) \otimes \qit/\zit')^{G_{k}} \arrow{r} & H^{2m}_{\text{\'et}} (\overline{X}, \qit/\zit'(m))^{G_{k}} \arrow{r} & Br^{m} (\overline{X})^{G_{k}}[\frac{1}{p}] \arrow[two heads]{r} & K, 
\end{tikzcd}
\]
where $K:= \text{ker } \{H^{1} (k, A^{m} (\overline{X})\otimes \qit/\zit') \to H^{1} (k, H^{2m}_{\text{\'et}} (\overline{X}, \qit/\zit'(m)))\}$. Since the middle vertical arrow has finite cokernel by Corollary \ref{natural}, what remains is to show that $K$ is finite. Now, $K$ is a subquotient of $Br^{m} (\overline{X})[\frac{1}{p}]$, which is of cofinite type. By Lemma \ref{nat-lem}, it thus suffices to prove that $K$ has finite exponent. However, by Corollary \ref{splitting}, the map $A^{m} (\overline{X}) \otimes \qit/\zit' \hookrightarrow H^{2m}_{\text{\'et}} (\overline{X}, \qit/\zit'(m))$ is split-injective in the isogeny category of $G_{k}$-modules. Hence, applying the functor $H^{1} (k, -)$ gives a split-injective map in the isogeny category of Abelian groups. Lemma \ref{lem-bas} then shows that $K$ is of finite exponent, as desired.

\Addresses
\end{document}